\documentclass[12pt]{amsart}
\usepackage{amssymb}
\newtheorem{theorem}{Theorem}

\newtheorem{corollary}[theorem]{Corollary}

\newtheorem{example}[theorem]{Example}

\newtheorem{lemma}[theorem]{Lemma}

\newtheorem{proposition}[theorem]{Proposition}
\newtheorem{remark}[theorem]{Remark}

\def\Label#1{\label{#1}}
\def \bR{\mathbb R}
\def \bC{\mathbb C}

\begin{document}

\title[The Goursat problem for a generalized Helmholz operator]
{The Goursat problem for a generalized Helmholz operator in the
plane}
\author[Peter Ebenfelt and Hermann Render]{Peter Ebenfelt and Hermann
Render}


\thanks{The first author is
supported in part by DMS-0401215.
The second author is
supported in part by Grant MTM2006-13000-C03-03 of the D.G.I. of Spain.}

\address{P. Ebenfelt: Department of Mathematics, University of California,
San Diego, La Jolla, CA 92093--0112, USA.}
\email{pebenfelt@math.ucsd.edu}
\address{H. Render: Departamento de Matem\'{a}ticas y Computaci\'{o}n,
Universidad de La Rioja, Edificio Vives, Luis de Ulloa s/n., 26004
Logro\~{n}o, Espa\~{n}a.} \email{render@gmx.de}

\maketitle

\begin{abstract} We consider the Goursat problem in the plane for partial
differential operators whose principal part is the $p$th power of
the standard Laplace operator. The data is posed on a union of $2p$
distinct lines through the origin. We show that the solvability of
this Goursat problem depends on Diophantine properties of the
geometry of lines on which the data is posed.
\end{abstract}

\section{Introduction} Let us consider in $\bR^2$ the
mixed Cauchy problem
\begin{equation}\Label{cauchy0}
\left\{
\begin{aligned}
&\Delta^p u + \sum_{|\alpha|\leq k_0} a_\alpha
\frac{\partial^{|\alpha|}u}{\partial x^\alpha}=f\\& P|(u-g),
\end{aligned}
\right.
\end{equation}
where $p$ is a positive integer, $k_0$ is an integer with $0\leq
k_0\leq 2p-1$, $\Delta$ denotes the standard Laplace operator in
$\bR^2$
$$\Delta:=\frac{\partial^2}{\partial x^2}+\frac{\partial^2}{\partial
y^2},$$ the coefficients $a_{\alpha}=a_{\alpha}(x,y)$ as well as the
data functions $f=f(x,y)$ and $g=g(x,y)$ are real-analytic functions
near $0$, and $P=P(x,y)$ is a homogeneous polynomial of degree $2p$.
Here, the notation $P|(u-g)$ means that $P$ divides $u-g$ in the
ring of germs of real-analytic functions at $0$. For instance, if
$P(x,y)=L(x,y)^{2p}$ for some linear function $L(x,y)$ (which is
equivalent to saying that the zero set of $P(x,y)$ consists of a
single line with multiplicity $2p$), then \eqref{cauchy0} with
$k_0=2p-1$ is a standard Cauchy problem with data on the line
$\{L(x,y)=0\}$ and the classical Cauchy-Kowalevsky Theorem
guarantees that \eqref{cauchy0} has a unique real-analytic solution
$u$ near $0$ for every choice of data functions $f$ and $g$. In the
recent paper \cite{EbRe06}, the authors show that if $P$ is elliptic
(i.e.\ the zero set of $P(x,y)$ consists of only the origin), then
\eqref{cauchy0} with $k_0=p$ has a unique solution $u$ for every
choice of data functions $f$ and $g$. In this paper, we shall
consider the case where the zero set of $P(x,y)$ is a union of $2p$
distinct lines (in which case \eqref{cauchy0} may be called a
Goursat problem). This case is much more subtle and leads to a small
divisor problem. We shall give a sufficient condition (which is also
necessary in the case $p=1$; see Section \ref{nec}) on the divisor
$P$ (see Theorem \ref{homodelp} below) for the homogeneous Goursat
problem
\begin{equation}\Label{goursatp1}
\left\{
\begin{aligned}
&\Delta^p u=f\\& P|(u-g)
\end{aligned}
\right.
\end{equation}
to have a unique real-analytic solution $u$ for every real-analytic
data $f$ and $g$. We shall also give a sufficient condition on $P$
(Theorem \ref{helmdelp} below) for the perturbed Goursat problem
\begin{equation}\Label{goursatp2}
\left\{
\begin{aligned}
&\Delta^p u+cu=f\\& P|(u-g),
\end{aligned}
\right.
\end{equation}
where $c=c(x,y)$ is a real-analytic function near $0$, to have a
unique real-analytic solution $u$ for every data function $f$ and
$g$.

The conditions on $P$ in Theorems \ref{homodelp} and \ref{helmdelp}
involve Diophantine properties of a determinant constructed from the
geometry of the lines constituting the zero set of $P$. For
instance, if $p=1$, so that $P$ has degree two and its zero set
consists of two distinct lines, then the condition can be phrased in
terms of the (acute) angle $\theta=2\pi \alpha$ between the two
lines. The necessary and sufficient condition for the homogeneous
Goursat problem
\begin{equation}\Label{goursat11}
\left\{
\begin{aligned}
&\Delta u=f\\& P|(u-g)
\end{aligned}
\right.
\end{equation}
to be solvable (Corollary \ref{homodel1}) is that
\begin{equation}\Label{dioleray}
\liminf_{\mathbb Z\ni m\to\infty} \left(\inf_{n\in \mathbb
Z}\left|\alpha-\frac{n}{m}\right|\right)>0,
\end{equation}
a condition that is satisfied by e.g.\ all non-Liouville numbers.
Our condition for the perturbed Goursat problem
\begin{equation}\Label{goursat12}
\left\{
\begin{aligned}
&\Delta u+cu=f\\& P|(u-g),
\end{aligned}
\right.
\end{equation}
to be solvable (Corollary \ref{helmdel1}) is more restrictive,
namely there exists a constant $C>0$ such that
\begin{equation}\Label{diophantine-1}
\left|\alpha-\frac{n}{m}\right|\geq\frac{C}{m^2} ,\quad
\forall n,m\in \mathbb Z, m\neq 0 .
\end{equation}
We note that every irrational number $\alpha$ that satisfies an
integral quadratic equation (like $\sqrt{k/l}$ for any integers $k$
and $l$) satisfies \eqref{diophantine-1} (by Liouville's Theorem on
Diophantine approximation).  We also point out that every
irrational, algebraic number satisfies
\begin{equation}\Label{diophantine-2}
\left |\alpha-\frac{n}{m}\right |\geq \frac{C_\mu}{m^\mu},\quad
\forall n,m\in \mathbb Z, m\neq 0 ,
\end{equation}
for some constant $C_\mu$ (that depends on $\mu$) and {\it every}
$\mu>2$ by the Thue-Siegel-Roth Theorem \cite{Ro55}). However, there
are algebraic numbers that do not satisfy \eqref{diophantine-1}.

We also mention that it follows from our proof that
\eqref{goursat12} has a unique formal power series solution for all
$f$ and $g$ if and only if $\alpha$ is irrational. Thus, as a
consequence of our results, we conclude that the family of Goursat
problems \eqref{goursat12}, parametrized by the angle $2\pi\alpha$
between the two lines in the zero set of $P$, displays  "chaotic"
behavior in that the set of parameters for which \eqref{goursat12}
is solvable is dense as is the set of parameters for which there is
not even a formal solution.

The homogeneous Goursat problem \eqref{goursat11} (i.e.\
\eqref{goursatp1} with $p=1$) can be transformed, by a simple linear
change of coordinates, into a Goursat problem considered by J. Leray
in \cite{L74}. (It was also briefly considered in its present form
by H. Shapiro in \cite{Shap89}.) Leray's main result is equivalent
to our Corollary \ref{homodel1}. The relationship between the two
Goursat problems and Leray's work is briefly explained in Section
\ref{lerayequiv} below. Leray's work was extended to complex
parameters and to higher dimensions by Yoshino in \cite{Y81a} and
\cite{Y81}. Other related work on mixed Cauchy and Goursat problems
include that of G\aa rding \cite{G65} (see also Theorem 9.4.2 in
H\"ormander \cite{H90}), Shapiro \cite{Shap89}, the first author and
Shapiro \cite{EbSh96}, and the authors \cite{EbRe06}. Our approach
to studying the Goursat problem is inspired by ideas from
\cite{Shap89} and \cite{EbSh96}. The proofs are based on a new
estimate for an associated Fischer operator in the real Fischer norm
(Theorem \ref{estimate}). The real Fischer norm was introduced in
\cite{Rend05} and was also used in \cite{EbRe06}.

This paper is organized as follows. We present our main results in
Section \ref{mainresults}. In Section \ref{lerayequiv}, we discuss
the relation between our results in the case $p=1$ and $c\equiv 0$
and those of Leray in \cite{L74}. An associated Fischer operator,
which is used in the proofs of the main results, is introduced in
Section \ref{s:est} and a crucial estimate is proved for that
operator (Theorem \ref{estimate}). The proof of Theorem
\ref{homodelp} is also given in that section. The proof of Theorem
\ref{helmdelp} is given in the subsequent section. In Section
\ref{ex}, we consider the case $p=2$ and present an explicit family
of examples to which Theorem \ref{helmdelp} can be applied (see
Theorem \ref{helmdel2}). Finally, in Section \ref{nec}, we show that
our condition in Corollary \ref{homodel1} is also necessary in this
case ($p=1$).

\section{Main results}\Label{mainresults}

We shall now formulate our results more precisely. We must first
introduce some notation. Let $B_{R}:=\left\{ (x,y)\in
\mathbb{R}^{2}:x^2 +y^2 <R^2 \right\} $ be the open disk of
radius $R$ in $\mathbb{R}^{2}$ (where $0<R\leq \infty ).$ We denote
by  $A\left( B_{R}\right) $ the algebra of all infinitely
differentiable functions $f:B_{R}\rightarrow \mathbb{C}$ such that
for any compact subset $ K\subset B_{R}$ the homogeneous Taylor
series $\sum_{m=0}^{\infty }f_{m}\left( x,y\right) $ converges
absolutely and uniformly to $f$ on $K$; here, $f_{m}$ is the
homogeneous polynomial of degree $m$ defined by the Taylor series of
$f$
\begin{equation*}
f_{m}\left( x,y\right) =\sum_{k+l =m}\frac{1}{k!
l !}\frac{%
\partial ^{m }f}{\partial x^{k}\partial y^l}\left( 0\right) x^{k }y^l.
\end{equation*}
Note that the functions in $A(B_R)$ are real-analytic. For a real
number $a$, we shall define the unimodular complex number
\begin{equation}\Label{A}
A=A(a):=\frac{a+i}{a-i}.
\end{equation}
As $a$ goes from $-\infty$ to $\infty$, $A$ ranges over the unit
circle (from 1 to 1 in the negative direction) and, hence, there is
a unique $\beta\in (0,1)$ such that $A=e^{2\pi i \beta}$. Note that
$\beta$ is rational precisely when $A$ is a root of unity. For
future reference, we observe, using the fact that $2\arctan
a=i\log(1-ia)/(1+ia)$, that for $a\in [0,\infty)$ the acute angle
between the lines $y=0$ and $x-ay=0$ is $\pi \beta$. Now, let us fix
a positive integer $p$, distinct real numbers $a_1, a_2, \ldots,
a_{2p-1}$, and write $a$ for the vector $a=(a_1,\ldots,a_{2p-1})$.
We shall denote by $P_a(x,y)$ the divisor
\begin{equation}\Label{Pa}
P_{a}(x,y):=y\, \Pi_{j=1}^{2p-1}(x-a_jy).
\end{equation}
If the divisor $P$ in \eqref{cauchy0} is a homogeneous polynomial of
degree $2p$ with $2p$ distinct lines as its zero set, then there is
no loss of generality in assuming that  $P$ is of the form
\eqref{Pa}, since the Laplace operator is rotationally invariant. We
associate to the vector $a$ a sequence of $2p\times 2p$ matrices
$\{M_{m,p,a}\}_{m=0}^\infty$, where
\begin{equation}\Label{Ma}
M_{m,p,a}:=\begin{pmatrix} 1& 1& 1 & \ldots & 1 &
1&\ldots& 1\\
1& A_1& A_1^2 & \ldots & A_1^{p-1} &
A_1^{m+p+1}&\ldots& A_1^{m+2p}\\
1& A_2& A_2^2 & \ldots & A_2^{p-1} & A_2^{m+p+1}&\ldots&
A_2^{m+2p} \\
\vdots & \vdots & \vdots & \ddots &\vdots & \vdots & \ddots &\vdots \\
1 & A_{2p-1} & A_{2p-1}^2 & \ldots & A_{2p-1}^{p-1} &
A_{2p-1}^{m+p+1}&\ldots& A_{2p-1}^{m+2p}
\end{pmatrix}.
\end{equation}
Here, $A_j:=A(a_j)$ where $A(a_j)$ is given by \eqref{A}. We shall
consider the Goursat problem
\begin{equation}\Label{goursatp}
\left\{
\begin{aligned}
& \Delta^pu+cu=f\\& P_a|(u-g),
\end{aligned}
\right.
\end{equation}
where the coefficient $c=c(x,y)$ as well as the data functions
$f=f(x,y)$, $g=g(x,y)$  belong to $A(B_R)$. Our first result
concerns the homogenenous problem, i.e. $c\equiv 0$.

\begin{theorem}\Label{homodelp}
Let $p$ be a positive integer and $a_1,\ldots,a_{2p-1}$
 real, distinct, non-zero numbers. Let $A_j:=A(a_j)$, for $j=1,\ldots, 2p-1$,
 be the unimodular complex numbers given by \eqref{A}, $P_a(x,y)$
 the homogeneous polynomial given by \eqref{Pa}, and $\{M_{m,p,a}\}_{m=0}^\infty$
 given by \eqref{Ma}. If $\det M_{m,p,a}\neq
 0$ for all integers $m\geq 0$, and
 \begin{equation}\Label{leraycond1}
\tau:=\liminf_{m\to\infty}\left(\det M_{m,p,a}\right )^{1/m}>0,
 \end{equation}
then the homogeneous Goursat problem
\begin{equation}\Label{goursatp0}
\left\{
\begin{aligned}
& \Delta^pu=f\\& P_a|(u-g)
\end{aligned}
\right.
\end{equation}
has a unique solution $u\in A(B_{\tau R})$ for every $f,g\in
A(B_R)$.
\end{theorem}

\begin{remark}\Label{rmkmatrix}{\rm For future reference, we note the following
identity
\begin{multline}\Label{Ma3}
\det M_{m,p,a}=\\\det \begin{pmatrix} A_1-1 & A_1^2-1 & \ldots &
A_1^{p-1}-1 &
A_1^{m+p+1}-1&\ldots& A_1^{m+2p}-1\\
A_2-1 & A_2^2-1 & \ldots & A_2^{p-1}-1 & A_2^{m+p+1}-1&\ldots&
A_2^{m+2p}-1\\
\vdots & \vdots & \ddots &\vdots & \vdots & \ddots &\vdots \\
 A_{2p-1}-1 & A_{2p-1}^2-1 & \ldots & A_{2p-1}^{p-1}-1 &
A_{2p-1}^{m+p+1}-1&\ldots& A_{2p-1}^{m+2p}-1
\end{pmatrix}.
\end{multline}
In particular, for $p=1$, we have $\det M_{m,p,a}=A_1^{m+2}-1$. }
\end{remark}

We mention that e.g.\ all numbers $a_1,\ldots,a_{2p-1}$ such that
$A_1,\ldots, A_{2p-1}$ are algebraic and $\det M_{m,p,a}\neq 0$ for
all $m$ satisfy \eqref{leraycond1} (see \cite{W00}, Lemma 2.1).

It will follow from our proof of Theorem \ref{helmdelp} below that
the Goursat problem \eqref{goursatp}, and hence in particular
\eqref{goursatp0}, has a unique formal solution $u$ if and only if
$\det M_{m,p,a}\neq 0$ for all integers $m\geq 0$. The Diophantine
condition \eqref{leraycond1} is sufficient (and necessary for $p=1$;
see Section \ref{nec} below) for the formal solution to
\eqref{goursatp0} to converge. For the formal solution to the
general Goursat problem \eqref{helmdelp} to converge, we need a
stronger condition. We have the following result.

\begin{theorem}\Label{helmdelp}
Let $p$ be a positive integer and $a_1,\ldots,a_{2p-1}$
 real, distinct, non-zero numbers. Let $A_j:=A(a_j)$, for $j=1,\ldots, 2p-1$,
 be the unimodular complex numbers given by \eqref{A}, $P_a(x,y)$
 the homogeneous polynomial given by \eqref{Pa}, and $\{M_{m,p,a}\}_{m=0}^\infty$
 given by \eqref{Ma}. If there exists a constant $C>0$
 such that
 \begin{equation}\Label{leraycond2}
\det M_{m,p,a}\geq \frac{C}{m^p},
 \end{equation}
for all natural numbers  $m \geq 1 $ then there exists $0<r\leq R$ such that the Goursat problem
\eqref{goursatp} has a unique solution $u\in A(B_{r})$ for every
$f,g\in A(B_R)$.
\end{theorem}

In Section \ref{ex} below, we give some explicit examples of
$a_1,a_2,a_3$ such that \eqref{leraycond2} holds for the
corresponding unimodular numbers $A_1,A_2,A_3$.

In the case $p=1$, the zero set of $P_a$ is the union of the two
distinct lines given by $y=0$ and $x=ay$. By the rotational symmetry
of $\Delta$, we may also assume that $a\geq 0$. If we denote the
acute angle between the two lines by $2\pi \alpha$ and by $\beta\in
(0,1/2]$ the number such that $A:=A(a)=e^{2\pi i\beta}$, then as
mentioned in the beginning of this section we have $\beta=2\alpha$.
As noted in Remark \ref{rmkmatrix} above, we have $\det
M_{m,p,a}=A^{m+2}-1$. The condition $\det M_{m,p,a}=A^{m+2}-1\neq 0$ is clearly
equivalent to $\alpha$ being irrational. Since
$$
|A^{m+2}-1|\approx\inf_{n\in \mathbb Z} |2\pi (m+2) \beta-2\pi n|=2\pi (m+2)
\inf_{n\in \mathbb Z} \left|\beta- \frac{n}{m+2}\right|,
$$
where by $E_k\approx F_k$ we mean $CF_k\leq E_k\leq DF_k$ for
nonzero constants $C,D$, it is not difficult to see that Theorems
\ref{homodelp} and \ref{helmdelp}, specialized to the case $p=1$,
can be formulated as follows.

\begin{corollary}\Label{homodel1} Let $\Gamma_1$, $\Gamma_2$ be two distinct lines
through the origin in $\bR^2$, and denote by $\theta=2\pi\alpha$ the
acute angle between them. Suppose that $\alpha$ is irrational and
satisfies the condition
\begin{equation}\Label{leraycond3}
\tau:=\liminf_{m\to\infty}\left(\inf_{n\in\mathbb Z}\left
|\alpha-\frac{n}{m}\right |\right)^{1/m}>0.
\end{equation}
Then, the homogeneous Goursat problem
\begin{equation}\Label{goursat10}
\left\{
\begin{aligned}
& \Delta u=f\\& u=g \quad \text{{\rm on $\Gamma_1\cup\Gamma_2$}}
\end{aligned}
\right.
\end{equation}
has a unique solution $u\in A(B_{\tau R})$ for every $f,g\in
A(B_R)$.
\end{corollary}

The condition \eqref{leraycond3} is also necessary for the
conclusion of Corollary \ref{homodel1} to hold. This fact is proved
in Section \ref{nec} below.  As mentioned in the introduction,
Corollary \ref{homodel1} is equivalent to the result of Leray in
\cite{L74}. A more detailed explanation of this equivalence is given
in Section \ref{lerayequiv} below.

We conclude this section by reformulating Theorem \ref{helmdelp} in
the case $p=1$.

\begin{corollary}\Label{helmdel1} Let $\Gamma_1$, $\Gamma_2$ be two distinct lines
through the origin in $\bR^2$, and denote by $\theta=2\pi\alpha$ the
acute angle between them. Suppose that $\alpha$ satisfies the
Diophantine condition
\begin{equation}\Label{diophantine}
\left |\alpha-\frac{n}{m}\right |\geq \frac{C}{m^2},\quad \forall
n,m\in \mathbb Z, m\neq 0
\end{equation}
for some constant $C>0$. Then, for any $c\in A(B_R)$, there exists
$0<r\leq R$ such that the Goursat problem
\begin{equation}\Label{goursatp10}
\left\{
\begin{aligned}
& \Delta u+cu=f\\& u=g \quad \text{{\rm on $\Gamma_1\cup\Gamma_2$}},
\end{aligned}
\right.
\end{equation}
 has a unique solution $u\in A(B_r)$ for every
$f,g\in A(B_R)$.
\end{corollary}

\section{Leray's Goursat problem}\Label{lerayequiv}

Consider the homogeneous Goursat problem
\begin{equation}\Label{goursat0}
\left\{
\begin{aligned}
&\lambda \frac{\partial^2u}{\partial x\partial y}+\Delta u=f\\&
xy|(u-g),
\end{aligned}
\right.
\end{equation}
where $\lambda$ is a real constant. It follows from the general
theory of Goursat (or mixed Cauchy) problems that \eqref{goursat0}
has a unique real-analytic solution near $0$, for all $f$ and $g$,
if $|\lambda|>2$ (see G\aa rding \cite{G65}; see also Theorem 9.4.2
in H\"ormander \cite{H90}). The case where $\lambda\in[-2,2]$ is
much more subtle, and was analyzed by Leray in \cite{L74} (see also
the work of Yoshino \cite{Y81a}, \cite{Y81} for extensions to
complex parameters and higher dimensions). For $\lambda\in [-2,2]$,
let $\beta\in [-1/4,1/4]$ denote the angle such that
$\lambda=2\sin(2\pi\beta)$. Leray showed that the unique solvability
of \eqref{goursat0} depends on Diophantine properties of $\beta$.
For instance, there is a unique formal power series solution $u$ for
every $f$ and $g$ if and only if $\beta$ is irrational. Leray also
gave a necessary and sufficient Diophantine condition on irrational
$\beta$ quaranteeing that this formal solution $u$ converges for all
convergent $f$ and $g$,
\begin{equation}\Label{dioleray2}
\liminf_{\mathbb Z\ni m\to\infty} \left(\inf_{n\in \mathbb
Z}\left|\beta-\frac{n}{m}\right|^{1/m}\right)>0.
\end{equation}
Let us show that this result, for $\lambda\in(-2,2)$, is equivalent
to our Corollary \ref{homodel1} above. Consider the linear change of
variables
\begin{equation}\Label{trans}
x \to - \sqrt{1-\frac{\lambda^2}{4}}x+\frac{\lambda}{2}y.
\end{equation}
As the reader can easily verify, this change of variables leads to
the following transformation for the principal symbol of the
operator
\begin{equation}
\lambda \frac{\partial^2}{\partial x\partial y}+\Delta  \to \Delta.
\end{equation}
Hence, the Goursat problem \eqref{goursat0} is transformed into the
following
\begin{equation}\Label{goursat01}
\left\{
\begin{aligned}
&\Delta u=f\\& y(x-ay)|(u-g),
\end{aligned}
\right.
\end{equation}
where
\begin{equation}\Label{b}
a:=\frac{\lambda/2}{\sqrt{1-(\lambda/2)^2}}.
\end{equation}
If we let $\theta=2\pi\alpha$ denote the acute angle between the two
lines $L_1:=\{y=0\}$ and $L_2:=\{x=by\}$ and $\beta$ the angle such
that $\lambda:=2\sin(2\pi\beta)$, then we have
$$\alpha=\frac{1-2\beta}{4}.
$$
Clearly, we have
$$
\liminf_{\mathbb Z\ni m\to\infty} \left(\inf_{n\in \mathbb
Z}\left|\beta-\frac{n}{m}\right|\right)^{1/m}=\liminf_{\mathbb Z\ni
 m\to\infty} \left(\inf_{n\in \mathbb
Z}\left|\alpha-\frac{n}{m}\right|\right)^{1/m}.
$$
This shows, as mentioned in the introduction, that Leray's result,
with $\lambda\in (-2,2)$, is equivalent to our Corollary
\ref{homodel1}, with $0<a<\infty$.

\section{An estimate for an associated Fischer
operator and the proof of Theorem $\ref{homodelp}$}\Label{s:est}

Let $\bC[x,y]$ denote the space of polynomials in $x,y$ with complex
coefficients. For each integer $m\geq 0$, we shall let $\mathcal
P_m$ denote the subspace of homogeneous polynomials of degree $m$.
We endow $\bC[x,y]$ with the real Fischer inner product
\begin{equation}
\langle f,g\rangle :=\int_{\bR^2}
f(x,y)\overline{g(x,y)}e^{-(x^2+y^2)}dxdy,
\end{equation}
and denote by $\|\cdot\|$ the corresponding norm (see
\cite{Rend05}). We shall fix a positive integer $p$ and distinct
real numbers $a_1,\ldots, a_{2p-1}$ and consider the Fischer
operator $F_a(q):=\Delta^p(P_aq)$, where $P_a$ is given by
\eqref{Pa}. Observe that $F_a$ is a linear operator sending
$\mathcal P_m$ into $\mathcal P_m$. Our main result in this section
is the following, in which the notation introduced above is used.

\begin{theorem}\Label{estimate} Let $p$ be a positive integer and $a_1,\ldots,a_{2p-1}$
 real, distinct, non-zero numbers. Let $A_j:=A(a_j)$, for $j=1,\ldots, 2p-1$,
 be the unimodular complex numbers given by \eqref{A} and $P_a(x,y)$
 the homogeneous polynomial given by \eqref{Pa}. Then the Fischer operator
 $F_a\colon \mathcal P_m\to\mathcal P_m$, for $m\geq 0$, is a
bijection if and only if $\det M_{m,p,a}\neq 0$, where $M_{m,p,a}$ is
given by \eqref{Ma}. Moreover, if $\det
M_{m,p,a}\neq 0$, then we have the estimate
\begin{equation}
\|P_aq\|\leq \frac{C} {|\det M_{m,p,a}|}\|\Delta^p(P_aq)\|,\quad
\forall q\in \mathcal P_m,
\end{equation}
for some $C\geq 0$ (independent of $m$).
\end{theorem}

For the proof of Theorem \ref{estimate}, we shall need the following
lemma. To state the lemma, we observe the well known fact that any
homogeneous polynomial $f(x,y)$ of degree $m$ can be expressed in
the following way
\begin{equation}\Label{expaninz}
f(x,y)=\sum_{k+l=m}f_{kl}z^k\bar z^l,
\end{equation}
where $z=x+iy$ and $\bar z=x-iy$.

\begin{lemma}\Label{Fischer} Let $f(x,y)$ be a homogeneous
polynomial of degree $m$ given by \eqref{expaninz}. Then, we have
\begin{equation}\Label{normid}
\|f\|^2=\pi m!\sum_{k+l=m}|f_{kl}|^2.
\end{equation}
\end{lemma}

\begin{proof} As in \cite{Rend05} (see also \cite{EbRe06}), we observe
that for any homogeneous polynomial $f(x,y)$ of degree $m$, we have
\begin{equation}\Label{norm1}
\|f\|^2=I_{2m+1}\int_{\mathbb T}|f(\eta)|^2ds_{\eta}
\end{equation}
where $\mathbb T$ denotes the unit circle in $\bR^2$, $ds$
arclength, and $I_{k}$ the integral
\begin{equation*}
I_{k}:=\int_{0}^{\infty }e^{-r^{2}}r^{k}dr.
\end{equation*}
A simple substitution argument gives
\begin{equation}\Label{int1}
I_{2m+1}=\int_{0}^{\infty
 }e^{-r^{2}}r^{2m+1}dr=\frac{1}{2}\int_{0}^{\infty
}e^{-x}x^{m}dx=\frac{1}{2}m!.
\end{equation}
Substituting \eqref{expaninz} in \eqref{norm1}, using the
parametrization $z=e^{i\theta}$ for $\mathbb T$ and the identity
\eqref{int1}, yields
\begin{equation}\Label{norm2}
\|f\|^2=\frac{1}{2}m!\sum_{k+l=m}\sum_{i+j=m}f_{kl}\overline{f_{ij}}
\int_0^{2\pi}e^{i(k+j-l-i)\theta}d\theta,
\end{equation}
from which \eqref{normid} readily follows.
\end{proof}

\begin{proof}[Proof of Theorem $\ref{estimate}$] We fix $f\in
\mathcal P_m$ and consider the equation
\begin{equation}\Label{Faisf}
F_a(q):=\Delta^p(P_aq)=f,
\end{equation}
for $q\in \mathcal P_m$. Note that $q\in \mathcal P_m$ solves
\eqref{Faisf} if and only if $u=P_aq$ solves the Goursat problem
\begin{equation}\Label{goursat2}
\left\{
\begin{aligned}
& \Delta^p u=f\\ & u(x,0)=u(a_1y,y)\ldots u(a_{2p-1}y,y)=0.
\end{aligned}
\right. \end{equation} We shall look for $u$ of the form $u=v+w$,
where $w(x,y)=(x^2+y^2)^p s(x,y)$ for some $s\in \mathcal H_{m}$
such that
\begin{equation}\Label{classic}
\Delta^p w(x,y)=\Delta^p((x^2+y^2)^ps(x,y))=f(x,y)
\end{equation}
and $v\in \mathcal H_{m+2p}$ satifies
\begin{equation}\Label{goursat3}
\left\{
\begin{aligned}
& \Delta^p v=0\\ & v(x,0)=-w(x,0)\\ &v(a_jy,y)=-w(a_jy,y), \quad
j=1,\ldots 2p-1.
\end{aligned}
\right. \end{equation} It is well known that \eqref{classic} has a
unique solution $w(x,y)=(x^2+y^2)^ps(x,y)$ (see e.g.\ \cite{Shap89}
and references therein). Moreover, in view of the results in
\cite{EbRe06}, we have
\begin{equation}
\|w\|\leq C_1\|f\|
\end{equation}
for some constant $C_1>0$. Thus, to complete the proof of the
theorem it suffices to show that \eqref{goursat3} has a solution
$v\in \mathcal P_{m+2p}$ for every $f\in \mathcal P_m$ if and only
if $\det M_{m,p,a}\neq 0$, and that, in this case,
\begin{equation}\Label{goal1}
\|v\|\leq \frac{C}{|\det M_{m,p,a}|}\|f\|
\end{equation}
for some constant $C>0$. To this end, we shall actually need the
exact form of the solution to \eqref{classic}. Using $z=x+iy$ and
$\bar z=x-iy$, we may write \begin{equation}\Label{id1}
w(x,y)=W(z,\bar z)=z^p\bar z^p\sum_{k+l=m} s_{kl}z^k\bar
z^l=\sum_{k+l=m}s_{kl}z^{k+p}\bar z^{l+p}.
\end{equation}
We observe that $\Delta=4\partial^2/\partial z\partial \bar z$.
Thus, if we write $f(x,y)=\sum_{k+l=m} f_{kl}z^k\bar z^l$, then
\eqref{classic} is equivalent to
\begin{equation}\Label{id2}
s_{kl}=\frac{f_{kl}}{4^p(k+1)\ldots (k+p)(l+1)\ldots (l+p)},\quad
\forall\ k+l=m.
\end{equation}
Now, we note that every function $v(x,y)$ that satisfies
$\Delta^pv=0$ is of the form
\begin{equation}\Label{vform}
v(x,y)=\sum_{t=0}^{p-1}\left({\bar z}^t\phi_t(z)+z^t\psi_t(\bar
z)\right),
\end{equation} where $\phi_t(z)$ and $\psi_t(\bar z)$ are holomorphic
functions of $z$ and $\bar z$, respectively.  The function $v$ is a
homogeneous polynomial of degree $m+2p$ if and only if
$\phi_t(z)=b_{p-1-t}z^{m+2p-t}$ and $\psi_t(\bar z)=c_t \bar z^{m+2p-t}$,
for constants $b_{p-1-t}$ and $c_t$ and $t=0,\ldots, p-1$. Using
this notation, equation \eqref{goursat3}  is
equivalent to finding monomials
\begin{equation}\Label{id3}\phi_t(z)=b_{p-1-t}z^{m+2p-t},\quad \psi_t(\bar z)=c_t \bar
z^{m+2p-t}, \end{equation} for $t=0,1,\ldots,p-1$, such that
\begin{equation}\Label{get1}
\sum_{t=0}^{p-1}\left({x}^t\phi_t(x)+x^t\psi_t(x)\right)=-W(x,x)
\end{equation}
and
\begin{multline}\Label{get2}
\sum_{t=0}^{p-1}\left({((a_j-i)y)}^t\phi_t((a_j+i)y)+
((a_j+i)y)^t\psi_t((a_j-i)y)\right)=\\-W((a_j+i)y,(a_j-i)y),\quad
j=1,\ldots, 2p-1.
\end{multline}
In \eqref{get2}, we use the fact that $\phi_t $ is homogeneous of
degree $m+2p-t $ and divide the equation by $ (a_j -i )^{m+2p } $.
With $A_j:=A(a_j)$ and $A(a)$ given by \eqref{A}, the equation
becomes
\begin{equation}\Label{get3}
\sum_{t=0}^{p-1}\left(A_j^{m+2p-t}y^t \phi_t(y )+ A_j^ty^t
\psi_t(y)\right)=-W(A_j y,y),\quad j=1,\ldots, 2p-1.
\end{equation}
Substituting \eqref{id2} and \eqref{id3} in \eqref{get1} and
\eqref{get3}, we obtain the following system of linear equations for
the coefficients $b_0,\ldots, b_{p-1}, c_0,\ldots, c_{p-1}$
\begin{equation}\Label{system1}
\begin{aligned}
\sum_{t=0}^{p-1}\left(b_{p-1-t}+c_t\right) &=-\sum_{k+l=m}\frac{f_{kl}}
{4^p(k+1)\ldots (k+p)(l+1)\ldots (l+p)}\\
\sum_{t=0}^{p-1}\left(A^{m+2p-t}_1 b_{p-1-t}+ A_1^t c_j\right)
&=-\sum_{k+l=m}\frac{f_{kl}A_1^{k} } {4^p(k+1)\ldots
(k+p)(l+1)\ldots
(l+p)}\\
&\vdots\\
\sum_{t=0}^{p-1}\left(A^{m+2p-t}_{2p-1} b_{p-1-t}+ A_{2p-1}^t
c_j\right) &=-\sum_{k+l=m}\frac{f_{kl}A_{2p-1}^{k}} {4^p(k+1)\ldots
(k+p)(l+1)\ldots (l+p)}
\end{aligned}
\end{equation}
If we write $d$ for the column vector of coefficients
$d=(c_0,\ldots,c_{p-1},b_0,\ldots, b_{p-1})^t$ and $e$ for the
column vector whose $(j+1)$th component, $j=0,\ldots,2p-1$, is given
by
$$
-\sum_{k+l=m}\frac{f_{kl}A_j^{k}} {4^p(k+1)\ldots (k+p)(l+1)\ldots
(l+p)},
$$
where we let $A_0:=1$, then \eqref{system1} can be written
\begin{equation}\Label{matrixeq}
M_{m,p,a}d=e, \end{equation}
 where $M_{m,p,a}$ is given by \eqref{Ma}.
We conclude, as claimed above, that \eqref{goursat3} has a unique
solution $v\in\mathcal P_{m+2p}$ for every
$f\in \mathcal P_m$ if and only if $\det M_{m,p,a}\neq 0$.

Let us now suppose that $\det M_{m,p,a}\neq 0$ and write $d_i$ for
the $i$th component of $d$, $i=1,\ldots, 2p$. Using Cramer's rule
and the fact that $|A_j|=1$, we conclude from \eqref{matrixeq} that
\begin{equation}\Label{coeffest}
|d_i|\leq C_1|\det M_{m,p,a}|^{-1}\sum_{k+l=m}\frac{|f_{kl}|}
{(k+1)\ldots (k+p)(l+1)\ldots (l+p)}.
\end{equation}
By the Cauchy-Schwarz inequality, we conclude that
\begin{equation}\Label{coeffest2}
|d_i|\leq C_1|\det
M_{m,p,a}|^{-1}\left(\sum_{k+l=m}|f_{kl}|^2\right)^{1/2} S_m,
\end{equation}
where $S_m$ denotes the sum
\begin{equation}\Label{Sm}
S_m:=\left(\sum_{k+l=m}\frac{1}{(k+1)^2\ldots (k+p)^2(l+1)^2\ldots
(l+p)^2}\right)^{1/2}.
\end{equation}
By setting $l=m-k$, we obtain
\begin{equation}\Label{Smest}
\begin{aligned}
S_m^2= & \sum_{k=0}^m\left(\prod_{j=1}^p(k+j)^2(m-k+j)^2\right)^{-1}\\
\leq &\,
2\sum_{k=0}^{[m/2]+1}\left(\prod_{j=1}^p(k+j)^2(m-k+j)^2\right)^{-1}\\
= &\,
2m^{-2p}\sum_{k=0}^{[m/2]+1}\left(\prod_{j=1}^p(k+j)^2\left((1+(j-k)/m\right)^2\right)^{-1}
\end{aligned}
\end{equation}
Now, note that, for $j=1,\dots,p$ and $k=0,\ldots,[m/2]+1$, we have
$(j-k)/m\geq -3/4$ when $m\geq 2$ and, hence, $(1+(j-k)/m)^{-2}\leq
16$. Consequently, we have
\begin{equation}\Label{Smest2}
S_m^2\leq \frac{32}{m^{2p}}
\sum_{k=0}^{[m/2]+1}\left(\prod_{j=1}^p(k+j)^2\right)^{-1}\leq
\frac{32}{m^{2p}} \sum_{k=0}^{\infty}\frac{1}{(k+1)^{2p}}\leq
\frac{C_2}{m^{2p}},
\end{equation}
for some $C_2>0$ independent of $m$. Thus, by Lemma \ref{Fischer},
we obtain from \eqref{coeffest2} and \eqref{Smest2} the following
estimates for the functions $\tilde \phi_t(z,\bar z):=\bar
z^t\phi_t(z)$, where $\phi_t$ is given by \eqref{id3},
\begin{equation}\Label{estphipsi}
\begin{aligned}
\|\tilde \phi_t\|=&\, \sqrt{(m+2p)!}\,|b_{p-1-t}|\\
\leq&\, C_1C_2|\det M_{m,p,a}|^{-1}\sqrt{(m+1)\ldots(m+2p)}\|f\| m^{-p}\\
\leq&\, C_3 |\det M_{m,p,a}|^{-1}\|f\|.
\end{aligned}
\end{equation}
We obtain a similar estimate for $\tilde \psi_t(z,\bar
z):=z^t\psi_t(\bar z)$.  These estimates yield \eqref{goal1} since
$v$ is given by \eqref{vform}. This completes the proof of Theorem
\ref{estimate}.
\end{proof}

The arguments in the proof above also yield a proof of Theorem
\ref{homodelp}. We conclude this section by giving this proof.

\begin{proof}[Proof of Theorem $\ref{homodelp}$] It is well known that to prove Theorem
\ref{homodelp} it suffices to show that the equation
\begin{equation}\Label{PDE}
\Delta^p(Pq)=f
\end{equation}
has a unique solution $q\in A(B_{\tau R})$ for every $f\in A(B_R)$
(see e.g.\ \cite{EbRe06}). As in the proof of Theorem
\ref{estimate}, we shall look for the solution $u:=P_aq$ in the form
$u=v+w$, where $w(x,y)=(x^2+y^2)^ps(x,y)$ satisfies \eqref{classic}
and $v$ solves \eqref{goursat3}. It is well known that $w\in A(B_R)$
(see \cite{Shap89}; see also \cite{EbRe06}). Thus, to complete the
proof, it suffices to show that $v\in A(B_{\tau R})$. We expand $v$
as a series $v=\sum_mv_m$, where the $v_m$ are the homogeneous
Taylor polynomials of degree $m$ of $v$. Similarly, we expand
$w=\sum_m w_m$ and $f=\sum_m f_m$. By homogeneity, we observe that
the homogeneous polynomials $v_m$, $w_m$, $f_m$ satisfy
\eqref{goursat3} (with $v=v_m$, $w=w_m$, and $f=f_m$). The fact that
$v\in A(B_{\tau R})$ now follows easily from the definition
\eqref{leraycond1} of $\tau$, the form \eqref{vform} of $v$, and the
estimate \eqref{estphipsi}. The details are left to the reader.
\end{proof}

\section{Proof of Theorem $\ref{helmdelp}$}

\begin{proof}[Proof of Theorem $\ref{helmdelp}$] We fix
$a=(a_1,\ldots, a_{2p-1})$ as in the theorem. For brevity, we denote
$P_a$ simply by $P$. To prove Theorem \ref{helmdelp}, it suffices to
show that
 there is $0<r\leq R$ such that the
equation
\begin{equation}\Label{PDE1}
(\Delta^p+c)(Pq)=f
\end{equation}
has a unique solution $q\in A(B_r)$ for every $f\in A(B_R)$. We
shall look for the solution $u=Pq$ as a series
$u=\sum_{m}u_m=\sum_{m}Pq_{m-2p}$, where the $u_m$ are the
homogeneous Taylor polynomials of degree $m$ of $u$. To this end, we
expand, similarly, both $f$ and $c$ as Taylor series $f=\sum_m f_m$
and $c=\sum_m c_m$. The equation \eqref{PDE} then implies
\begin{equation}\Label{basic0}
\Delta^p(Pq_j)=f_j,\quad j=0,1,\ldots, 2p-1,
\end{equation}
and, for each $m\geq 2p$,
\begin{equation}\Label{basic1}
\Delta^p(Pq_m)=f_m-\sum_{k=0}^{m-2p}c_{m-k-2p}Pq_{k}.
\end{equation}
Since the Fischer operator $F=F_a$, given by $F(q)=\Delta^p(Pq)$,
 is bijective $\colon \mathcal P_m\to \mathcal P_m$ for every $m$
(by Theorem \ref{estimate}), we can solve, uniquely, \eqref{basic0}
and \eqref{basic1} inductively for $q_m$. This gives us a unique
formal power series solution $u=\sum_mu_m$ with $u_m=Pq_{m-2p}$. It
remains to prove that there is $r>0$ such that this series converges
to a function in $A(B_{r})$. For this, we observe that Theorem
\ref{estimate} and the assumption \eqref{leraycond2} implies the
following estimate
\begin{equation}\Label{basic2}
\|u_{m+2p}\|\leq Cm^{p}\|\Delta(Pq_m)\|\leq Cm^{p}\left(
\|f_m\|+\sum _{k=0}^{m-2p}\|c_{m-k-2p}u_{k+2p}\|\right)
\end{equation}
To prove that $u\in A(B_{r})$, we must show (see Proposition 16 in
\cite{EbRe06}) that for every $0<\rho<r$ there is a constant $B>0$
such that
\begin{equation}\Label{ind}
\|u_k\|\leq B\rho^{-k}\sqrt{k!}
\end{equation}
for every $k\geq 0$. Let us pick $\rho<\sigma<R$. In view of
Proposition 16 in \cite{EbRe06}, we may assume that there are
constants $D$ and $E$ such that
\begin{equation}\Label{assump}
 \max_{\theta \in \mathbb{T}
}\left| c_{k}\left( \theta \right) \right| \leq D\sigma^{-k},\quad
\left\| f_{k}\right\| \leq  E \rho^{-k} \sqrt{k!},
\end{equation}
for all $k\geq 0$. We shall prove \eqref{ind} by induction. Thus,
assume that \eqref{ind} holds for all $k\leq m+2p-1$. We shall prove
that \eqref{ind} holds also for $k=m+2p$, provided that $m$ is large
enough. By using \eqref{assump}, the induction hypothesis, and
Proposition 8 in \cite{EbRe06} (see also Proposition 7 in that
paper), we conclude from \eqref{basic2} the following estimate, for
some constant $F>0$,
\begin{equation}
\begin{aligned} \|u_{m+2p}\| \leq &Cm^{p}\left(E\rho^{-m}\sqrt{m!}+
\sum _{k=0}^{m-2p}F\sigma^{-(m-k-2p)}[(k+2p+1)\ldots (m-1)m]^{1/2}\|u_{k+2p}\|\right )\\
\leq &Cm^{\mu-1}\left(E\rho^{-m}\sqrt{m!}+ \sum
_{k=0}^{m-2p}BF\sigma^{-(m-k-2p)}\rho^{-(k+2p)}\sqrt{m!}\right )\\=&
B\rho^{-(m+2p)}\sqrt{(m+2p)!}\, T_m,
\end{aligned}
\end{equation}
where
\begin{equation}
\begin{aligned}
T_m: &=Cm^{p}\frac{\rho^{2p}}{\sqrt{(m+1)(m+2)}}\left
(E/B+F\sum_{k=0}^{m-2p} \left(\frac{\rho}{\sigma}\right)^{m-k-2p}\right)\\
&\leq Cm^{p}\frac{\rho^{2p}}{\sqrt{(m+1)\ldots (m+2p)}}\left
(E/B+F\frac{1}{1-\rho/\sigma}\right).
\end{aligned}
\end{equation}
Since $\rho<r$, we can make $T_m\leq 1$ for all $m$ by requiring
$0<r\leq R$ small enough (and keeping $\sigma<R$ fixed). This proves
Theorem \ref{helmdelp}.
\end{proof}

\section{Examples of solvable Goursat problems for
$\Delta^2+c$}\Label{ex}

In this section, we shall consider the following one-parameter
family of Goursat problems
\begin{equation}\Label{goursat22}
\left\{
\begin{aligned}
&\Delta^2 u+cu=f\\& P_t|(u-g),
\end{aligned}
\right.
\end{equation}
where $P_t(x,y)$, for $t> 0$, denotes the divisor
\begin{equation}\Label{Pt}
P_t(x,y):=xy(x-ty)(x-y/t).
\end{equation}
Recall that $A=A(t)$ denotes the unimodular number given by
\eqref{A} (with $a=t$). Let us denote by $\beta=\beta(t)$ the number
$\beta\in (0,2\pi)$ such that $A=e^{2\pi i \beta}$. We shall prove
the following result.

\begin{theorem}\Label{helmdel2} Let $t> 0$ and $\beta:=\beta(t)$ as defined above.
Suppose that
$\beta$ satisfies the Diophantine condition
\begin{equation}\Label{diophantine2}
\left |\beta-\frac{n}{m}\right |\geq \frac{C}{m^2},\quad \forall
n,m\in \mathbb Z, m\neq 0 ,
\end{equation}
for some constant $C>0$. Then, for any $c\in A(B_R)$, there exists
$0<r\leq R$ such that the Goursat problem \eqref{goursat22}
 has a unique solution $u\in A(B_r)$ for every
$f,g\in A(B_R)$.
\end{theorem}

Theorem \ref{helmdel2} is a direct consequence of Theorem
\ref{helmdelp}, with $p=2$, and the following proposition.

\begin{proposition}\Label{matrixcomp} Let $t>0$,
$a=(a_1,a_2,a_3):=(0,t,1/t)$,  and let $M_{m,p,a}$ be the matrix
defined by \eqref{Ma}  with $p=2$. If
$\beta=\beta(t)$ satisfies
\begin{equation}\Label{diophantine3}
\left |\beta-\frac{n}{m}\right |\geq \frac{C}{m^\mu},\quad \forall
n,m\in \mathbb Z, m\neq 0 ,
\end{equation}
for some constant $C>0$, then
\begin{equation}\Label{detAk}
\left|\det M_{m,p,a}\right| \geq \frac{D}{m^{2\mu -2}},
\end{equation}
for some $D>0$.
\end{proposition}

\begin{proof} It is easy to check that the unimodular numbers $(A_1,A_2,A_3)$ that correspond
to the vector $a$ is $(-1,A,B)$, where $AB=-1$ and, in view of the
discussion preceding Corollary \ref{homodel1},
\begin{equation}\Label{Am-1}
|A^m-1|\geq \frac{C'}{m^{\mu-1}}.
\end{equation}
(Of course, $A$ is given by \eqref{A}, but only the above two facts
will be needed in the proof.) To prove the proposition, it suffices,
in view of Remark \ref{rmkmatrix}, to show that $|N_m|\geq
C'/m^{2\mu-2}$, where
\begin{equation}\Label{newmatrix}
N_{m}:=M_{m-4,2,a} =\det \left(
\begin{array}{ccc}
-2 & (-1)^{m-1}-1 & (-1) ^{m}-1 \\
A -1 & A ^{m-1}-1 & A ^{m}-1 \\
B -1 & B ^{m-1}-1 & B ^{m}-1
\end{array}
\right) .
\end{equation}
We obtain, since $AB=-1$,
\begin{equation*}
A ^{m}N_{m}=\det \left(
\begin{array}{ccc}
-2 & (-1)^{m-1}-1 & (-1) ^{m}-1 \\
A -1 & A ^{m-1}-1 & A ^{m}-1 \\
-A ^{m-1}-A ^{m} & A (-1)^{m-1}-A ^{m} & (-1)^{m}-A ^{m}
\end{array}
\right) .
\end{equation*}
If $m$ is even, then
\begin{equation*}
A ^{m}N_{m}=\det \left(
\begin{array}{ccc}
-2 & -2 & 0 \\
A -1 & A ^{m-1}-1 & A ^{m}-1 \\
-A ^{m-1}-A ^{m} & -A -A ^{m} & 1-A ^{m}
\end{array}
\right) .
\end{equation*}
A straightforward calculation shows that
\begin{equation}\Label{even}
A^mN_M=4A(A^m-1)(A^{m-2}-1).
\end{equation}
If $m$ is odd, then
\begin{equation*}
A ^{m}N_{m}=\det \left(
\begin{array}{ccc}
-2 & 0 & -2 \\
A -1 & A ^{m-1}-1 & A ^{m}-1 \\
-A ^{m-1}-A ^{m} & A -A ^{m} & -1-A ^{m}
\end{array}
\right) .
\end{equation*}
This time we get
\begin{equation}\Label{odd}
A^mN_M=-2(A^{m-1}-1)^2(A^{2}+1).
\end{equation}
The conclusion $|N_m|\geq C'/m^{2\mu-2}$ follows easily from
\eqref{even} and \eqref{odd}. This completes the proof of the
proposition.
\end{proof}

\section{Divergence of formal solutions when $p=1$ and $\tau=0$.}\Label{nec}

 We now show that, for $p=1$ and irrational angles $\alpha$ between the two lines
$\Gamma_1$ and $\Gamma_2$, the formal solution $u$ to
\eqref{goursat10}, with $f$ convergent and $g\equiv 0$, need not
converge when $\tau$, given by \eqref{leraycond1}, is zero. A
discussion similar to the one that follows can also be found in
\cite{Shap89}. Using the notation and setup in the proof of Theorem
\ref{estimate}, let us choose $f$ such that for each $m$ we have,
for $k+l=m$,
\begin{equation}
f_{kl}=\left\{\begin{aligned}&R^{-m},\quad k=0\\&0,\quad k>0.
\end{aligned}\right.
\end{equation}
Note that $f\in A(B_R)$. Let us consider the Goursat problem
\eqref{goursat10} with $g=0$. By following the argument in the proof
of Theorem \ref{estimate} above, we conclude that the formal
solution is of the form $u=v+w$, where $w$ is the formal solution to
\eqref{classic} and $v(x,y)$ is the formal solution to
\eqref{goursat3}. Hence, $v$ is of the form
$v(x,y)=\phi(z)+\psi(\bar z)$. It is well known that the solution
$w$ to \eqref{classic} converges to a function in $A(B_R)$ (see
\cite{Shap89}; see also \cite{EbRe06}). Thus, the solution $u$ to
the Goursat problem converges if and only if the two power series
$\phi(z)=\sum_m b_{m}z^m$ and $\psi(\bar z)=\sum_m c_{m}\bar z^m$
converge. With $p=1$, it is easy to solve the system of equations
\eqref{system1} for $b_m$ and $c_m$ explicitly and we obtain
\begin{equation}
b_m=\frac{1}{(1-A^m)}\frac{A-1}{2R^{m-2}(m-1)},
\end{equation}
(A similar identity holds, of course, for $c_m$.) The radius of
convergence of the series $\phi(z)=\sum_m b_m z^m$ is
\begin{equation}
R\liminf_{m\to\infty} |1-A^m|^{1/m}=0,
\end{equation}
proving the assertion above that the solution $u$ does not converge.
We conclude this paper by giving an explicit example of a number
$\beta$ in $A=e^{2\pi i\beta}$ such that $\tau=0$.

\begin{example} {\rm Let us define
\begin{equation}
\beta:=\sum_{k=1}^\infty 10^{-p_k},
\end{equation}
where $p_k$ is defined recursively by $p_1=1$ and $p_{k+1}=p_k+k\,
10^{p_k}$. Note that, for every $N$, the rational number
$$r_N:=\sum_{k=1}^N10^{-p_k}=\frac{q_N}{10^{p_N}}$$ satisfies
$$
|\beta-r_N|\leq\frac{2}{10^{p_{N+1}}}.
$$
Consider the subsequence $m_N:=10^{p_N}$ and note that
$$
|A^{m_N}-1|\leq C \inf_{p,q\in \mathbb Z_+}q\left
|\beta-\frac{p}{q}\right|\leq
2\frac{10^{p_N}}{10^{p_{N+1}}}=\frac{2}{10^{p_{N+1}-p_N}}
$$
Thus, we have
$$
|A^{m_N}-1|^{1/m_N}\leq
\frac{C}{10^{(p_{N+1}-p_N)/10^{p_N}}}=\frac{C}{10^N}\to 0,
$$
which shows that $\tau=\liminf_{k\to\infty}|A^k-1|^{1/k}=0$. }
\end{example}

\end{document}